\documentclass[12pt]{article}
\usepackage{amsmath}
\usepackage{amssymb}
\begin{document}
\newtheorem{theorem}{Theorem}[section]
\newtheorem{remark}[theorem]{Remark}
\newtheorem{mtheorem}[theorem]{Main Theorem}
\newtheorem{observation}[theorem]{Observation}
\newtheorem{proposition}[theorem]{Proposition}
\newtheorem{lemma}[theorem]{Lemma}
\newtheorem{note}[theorem]{}
\newtheorem{1killlemma}[theorem]{First Killing-Lemma}
\newtheorem{2killlemma}[theorem]{Second Killing-Lemma}
\newtheorem{corollary}[theorem]{Corollary}
\newtheorem{notation}[theorem]{Notation}
\newtheorem{example}[theorem]{Example}
\newtheorem{definition}[theorem]{Definition}

\renewcommand{\labelenumi}{(\roman{enumi})}
\newcommand{\dach}[1]{\hat{\vphantom{#1}}}
\numberwithin{equation}{section}
\def\Z{{ \mathbb Z}}
\def\E{{ \mathbb E}}
\def\N{{ \mathbb N}}
\def\DD{ \mathbb D}
\def\GG{ \mathbb G}
\def\BBB{ B_\BB}
\def\D{{\hat{D}}}
\def\Q{{\mathbb Q}}
\def\G{\hat{G}}
\def\C{\hat{C}}
\def\T{{\cal T}}
\def\FF{{\mathfrak  F}}
\def\FFF{{\mathfrak  F}^*}
\def\C{{\mathfrak C}}
\def\X{{\mathfrak X}}
\def\Y{{\mathfrak Y}}
\def\RX{R\langle X \rangle}
\def\GX{\langle X \rangle}
\def\RXa{R\langle X_\alpha \rangle}
\def\RrXa{R(X_\alpha)}
\def\RrYa{R(Y_\alpha)}

\def\RaX{R_\alpha\langle X_\alpha \rangle}
\def\RraX{R_\alpha (X_\alpha)}
\def\K{{\mathfrak K}}
\def\qa{(R_\alpha, G_\alpha, \sigma_\alpha, \sigma_{\alpha *})}
\def\qb{(R_\beta, G_\beta, \sigma_\beta, \sigma_{\beta *})}
\def\R{\widehat{R}}
\def\F{\widehat{F}}
\def\G{\widehat{G}}
\def\T{{\cal T}}
\def\B{\widehat{B}}
\def\BC{\widehat{B_C}}
\def\BCC{\widehat{B_\C}}
\def\restr{\restriction}
\def\Aut{{\rm Aut\,}}
\def\Map{{\rm Map\,}}
\def\Im{{\rm Im\,}}
\def\ker{{\rm ker\,}}
\def\lg{{\rm lg\,}}
\def\br{{\rm br\,}}
\def\inf{{\rm inf\,}}
\def\sup{{\rm sup\,}}
\def\Br{{\rm Br\,}}
\def\Yphi{Y_{[\phi]}}
\def\Ypsi{Y_{[\psi]}}
\def\Xphi{X_{\tilde{\phi}}}
\def\Xpsi{X_{\tilde{\psi}}}
\def\a{\alpha}
\def\abar{\overline{\alpha}}
\def\aa{{\bf a}}
\def\to{\rightarrow}
\def\arr{\longrightarrow}
\def\sigmaa{{\bf \Sigma_a}}
\def\End{{\rm End\,}}
\def\Ines{{\rm Ines\,}}
\def\Hom{{\rm Hom\,}}
\def\Fin{{\rm Fin\,}}
\def\restr{\upharpoonright}
\def\Ext{{\rm Ext}\,}
\def\Hom{{\rm Hom}\,}
\def\End{{\rm End}\,}
\def\Aut{{\rm Aut}\,}
\def\ker{{\rm ker}\,}
\def\Ann{{\rm Ann}\,}
\def\defe{{\rm def}\,}
\def\rk{{\rm rk}\,}
\def\crk{{\rm crk}\,}
\def\nuc{{\rm nuc}\,}
\def\Dom{{\rm Dom}\,}
\def\Im{{\rm Im}\,}
\def\Yphi{Y_{[\phi]}}
\def\Ypsi{Y_{[\psi]}}
\def\Xphi{X_{\tilde{\phi}}}
\def\Xpsi{X_{\tilde{\psi}}}
\def\a{\alpha}
\def\abar{\overline{\alpha}}
\def\aa{{\bf a}}
\def\ra{\rightarrow}
\def\arr{\longrightarrow}
\def\iff{\Longleftrightarrow}
\def\sigmaa{{\bf \Sigma_a}}
\def\mm{{\mathfrak m}}
\def\X{{\mathfrak X}}
\def\Diam{\diamondsuit}

\title{{\sc Decompositions of Reflexive Modules}
\footnotetext{This work is supported by the project
No. G-545-173.06/97 of the German-Israeli
Foundation for Scientific Research \& Development\\
AMS subject classification:\\
primary: 13C05, 13C10, 13C13, 20K15,20K25,20K30  \\
secondary: 03E05, 03E35\\
Key words and phrases: almost free modules, reflexive modules,
duality theory, \\
modules with particular monomorphism \\ GbSh 716 in Shelah's list
of publications }}

\date{}
\author{ R\"udiger G\"obel and Saharon Shelah}
\maketitle
\begin{abstract} In a recent paper \cite{GS} we answered to the negative
a question raised in the book by Eklof and Mekler \cite [p. 455,
Problem 12]{EM} under the set theoretical hypothesis of $
\Diam_{\aleph_1}$ which holds in many models of set theory. The
Problem 12 in \cite{EM} reads as follows: {\it If $A$ is a dual
(abelian) group of infinite rank, is $A \cong A \oplus \Z$}?  The
set theoretic hypothesis we made is the axiom $ \Diam_{\aleph_1}$
which holds in particular in G\"odel's universe as shown by R.
Jensen, see \cite{EM}. Here we want to prove a stronger result
under the special continuum hypothesis (CH).

The question in \cite{EM} relates to dual abelian groups. We want
to find a particular example of a dual group, which will provide a
negative answer to the question in \cite{EM}. In order to derive a
stronger and also more general result we will concentrate on
reflexive modules over countable principal ideal domains $R$.
Following H. Bass \cite{B} an $R$-module $G$ is reflexive if the
evaluation map $\sigma: G \arr G^{**}$ is an isomorphism. Here
$G^* = \Hom (G, R)$ denotes the dual group of $G$. Guided by
classical results the question about the existence of a reflexive
$R$-module $G$ of infinite rank with $G \not\cong G \oplus R$ is
natural, see \cite [p. 455]{EM}. We will use a theory of bilinear
forms on free $R$-modules which strengthens our algebraic results
in \cite{GS}. Moreover we want to apply a model theoretic
combinatorial theorem from \cite{S} which allows us to avoid the
weak diamond principle. This has the great advantage that the
used prediction principle is still similar to $ \Diam_{\aleph_1}$
but holds under CH. This will simplify algebraic argument for $G
\not\cong G \oplus R$.
\end{abstract}

\section{Introduction}
Let $R$ be a countable principal ideal domain with $1 \ne 0$ and
$S$ a multiplicative closed subset of $R \setminus \{0\}$
containing $1$. If $S = \{ s_n : \ n \in \omega \}, \ s_0 = 1 $
and $q_n = \prod\limits_{i < n} s_i$, then

\begin{eqnarray} \label{qn-s}
q_{n+1} = q_n s_n \mbox{ and } \ \bigcap\limits_{n \in \omega}
q_nR = 0.
\end{eqnarray}

 The condition (\ref{qn-s}) requires that $R$ is not a
field, but then we may choose $S = R \setminus \{0\}$ or any other
`classical example'. The right ideals $q_n R$ form a basis of the
S-topology on $R$ which is Hausdorff by (\ref{qn-s}). The S-adic
completion of $R$ under the $S$-topology is the ring $\R$ which
has the size $2^{\aleph_0}$, see G\"obel, May \cite{GM} for
properties on $\R$. Similar consideration carry over to (free)
$R$-modules which we will use in the next sections. If $F$ is a
free $R$-module then the $S$-topology is generated by $F q_n \ (n
\in \omega)$ and if $\F$ denotes the S-adic completion of $F$ then
$F \subset \F$ and $F$ is pure and dense in $\F$. Recall that $F
\subseteq_* \F$ is pure (w.r.t. $S$) if and only if $$\F s \cap F
\subseteq F s \mbox{ for all } s \in S.$$ Also $F$ is dense in
$\F$ if and only if $\F / F$ is $S$-divisible (in the obvious
sense). Note that an element $e$ of an $R$-module $G$ is pure, we
write $e \in_* G$ if $e R \subseteq_* G$. If $X \leq G$, then $\GX
$ denotes the submodule generated by $X$ and $\GX_*$ denotes the
submodule purely generated by $X$. The following observation is
well-known and can be looked up in \cite{GM}.

\begin{observation} \label{obs1}
If $0 \neq r_n \in R$ infinitely often ($n \in \omega$), then we
can find $\epsilon_n \in \{0, 1 \}$ such that $$\sum\limits_{n \in
\omega} r_n q_n \epsilon_n \in \R \setminus R.$$
\end{observation}

 If $G$ is any $R$-module then $G^* = \Hom (G, R)$
denotes its dual module, and $G$ is a dual module if $G \cong D^*$
for some $R$-module $D$. Particular dual modules are reflexive
modules $D$ introduced by Bass \cite [p. 476]{B}. We need the
evaluation map $$\sigma = \sigma_D: D \arr D^{**} \ \ (d \arr
\sigma (d))$$ where $\sigma (d) \in D^{**}$ is defined by
evaluation $$ \sigma (d): D^* \arr R \ \ (\varphi \arr \varphi (d)
).$$ The module $D$ is {\it reflexive} if the evaluation map
$\sigma_D$ is an isomorphism. Obviously $D \cong (D^*)^*$ is a
dual module if $D$ is reflexive. However, there are very many
$\aleph_1$-free modules $G$ with $G^* = 0$, see \cite{CG}. An
$R$-module is $\aleph_1$-free if all countable submodules are
free. Using the special continuum hypothesis CH the existence of
many reflexive modules (Section 3) will follow from considerations
of free modules with bilinear form in Section 2. Classical
examples are due to Specker and {\L}os and it follows from
Specker's theorem that reflexive modules must be $\aleph_1$-free,
see Fuchs \cite{Fu}. Hence it is not too surprising that we will
deal with free $R$-modules first in Section 2. The bilinear form
is needed to control their duals when passing from $\aleph_0$ to
size $\aleph_1$ in Section 3.

In order to find reflexive groups $G$ of cardinality $\aleph_1${
with $G \not\cong R \oplus G$ we must discard all possible
monomorphisms
\begin{eqnarray} \label{discard}
\varphi : G \hookrightarrow G \mbox{ with } G \varphi \oplus eR =
G. \end{eqnarray} This will be established with the help of an
$\aleph_1$-filtration $$ G = \bigcup_{\alpha \in \omega_1}
G_\alpha \mbox{ of countable, pure, free submodules } G_\alpha$$
such that $G_\alpha$ is a summand with $R$-free complement of any
$G_\beta$ for $\alpha < \beta$ if $\alpha$ does not belong to a
fixed stationary subset $E$ of $\omega_1$. Given $\varphi$ as in
(\ref{discard}) by a back and forth argument there is a cub $C
\subseteq \omega_1$ such that $G_\alpha\varphi \subseteq G_\alpha$
and $e \in G_\alpha$ for all $\alpha \in C$.

If $\alpha \in E \cap C$ then
\begin{eqnarray}\label{discard-a} G_\alpha \varphi \oplus eR =
G_\alpha.\end{eqnarray} And if $\varphi\restr G_\alpha$ is
predicted by a function $\varphi_\alpha$ as under the assumption
of the diamond principle $ \Diam_{\aleph_1}$, then we construct
$G_{\alpha + 1}$ by a Step-Lemma from $G_{\alpha}$ and
$\varphi_\alpha$ such that $\varphi_\alpha$ does not extend to
$G_{\alpha + 1}$. Moreover $G_{\alpha + 1}$ is the $S$-adic
closure of $G_\alpha$
 in $G$ by the summand property mentioned above. Hence
$\varphi$ coincides with $\varphi_\alpha$ on $G_{\alpha}$ and maps
$G_{\alpha + 1}$ into itself, a contradiction. However assuming CH
only weaker prediction principles like weak diamond
$\Phi_{\aleph_1}$ are available, see  Devlin and Shelah \cite{DS}
or Eklof and Mekler \cite[p. 143, Lemma 1.7]{EM}. If we discard
$\varphi_\alpha$ in the construction as before, then
$\varphi\restr G_\alpha$ might extend because it does not entirely
agree with $\varphi_\alpha$. In this case $\varphi$ will show up
at some $\alpha < \beta \in E \cap C$ and we have a new chance to
discard $\varphi\restr G_{\alpha + 1}$ at level $\beta$ - provided
we know $\varphi\restr G_{\alpha + 1}$ when constructing $G_{\beta
+ 1}$.

This time we need a stronger algebraic algebraic Step-Lemma: Note
that $G_{\alpha + 1} \oplus F_\beta = G_\beta$. Now we have a
partial map $\varphi_\beta := \varphi\restr G_{\alpha + 1}$ with
domain $\Dom \varphi_\beta = G_{\alpha + 1} $ a summand of
$G_\beta$ and the splitting property (\ref{discard-a}) for the new
map $\varphi_\beta$:
\begin{eqnarray}\label{discard-b}
(G_{\alpha + 1}\varphi_\beta \oplus eR) \oplus F_\beta =
G_\beta.\end{eqnarray} In order to proceed we must discard these
partial maps and indeed we are able to prove a generalized
Step-Lemma taking care of partial maps (\ref{discard-b}) in
Section 2. Moreover, a counting argument also shows that a list of
such partial maps $$\{\varphi_\beta : \  \beta \in E \}$$ exists
which predict any given map $\varphi: G \arr G$ such that the
following holds.
\begin{lemma}\cite{S}\label{predict}
Let $ G = \bigcup_{\alpha \in \omega_1} G_\alpha$ be an
$\aleph_1$-filtration of $G$ and $E \subseteq \omega_1$ be a
stationary subset of $\omega_1$. Then there is a list of
predicting partial maps $$\{\varphi_\alpha :G_\alpha \arr G_\alpha
: \ \alpha \in E \}$$ such that for any countable subset $A$ of
$G$ there is an ordinal $\beta \in E$ (in fact an unbounded set of
such ordinals) with $$\varphi \restr A = \varphi_\beta\restr
(G_\beta \cap A).$$
\end{lemma}
For a suitable $A = G_{\alpha + 1}$ and $\varphi\restr A$ we
choose $\beta \in E$ by Lemma \ref{predict} with $\varphi\restr
G_{\alpha + 1} = \varphi_\beta \restr G_{\alpha + 1} $ and
$G_\beta = G_{\alpha + 1} \oplus F_\beta$. Then we construct
$G_{\beta + 1}$ from $G_\beta$ and the given $\varphi_\beta \restr
G_{\alpha + 1}$ such that $\varphi_\beta \restr G_{\alpha + 1}$
does not extend to $G_{\beta + 1}$. This contradiction will
provide the

\begin{mtheorem} \label{mainth}
(ZFC + CH) \  If $R$ is a countable domain but not a field, then
there is a family of $2^{\aleph_1}$ pair-wise non-isomorphic
reflexive $R$-modules $G$ of cardinality $\aleph_1$ such that $G
\not\cong R \oplus G$.
\end{mtheorem}

We also would like to draw attention to a slight modification of
the proof of the Main Theorem \ref{mainth}. In addition we may
assume that $G$ in the Main Theorem \ref{mainth} is essentially
indecomposable, that any decomposition into two summands has one
summand free of finite rank. This follows from a split realization
result Corollary \ref{ring} with $\End G \cong A\oplus \Fin(G)$
where $A$ is any $R$-algebra which is free of countable rank.
Recall that $\Fin(G)$ is the ideal of all endomorphisms of $G$ of
finite rank.

\bigskip

The formal proof for the prediction Lemma \ref{predict} is not
complicated and uses repeatedly often the weak diamond prediction
$\Phi_{\aleph_1}$. It is also clear from what we said that the
underlying module theory is not essential for proving Lemma
\ref{predict} and that it should be possible to replace modules by
many other categories like non-commutative groups, fields or
Boolean algebras. In order to cover all these possibilities, the
prediction principle is formulated in terms of model theory and
will appear in this setting in a forth coming book by Shelah
\cite[Chapter IX, Claim 1.5]{S}.

We close this introduction with some historical remarks. Using a
theorem of {\L}o{\`s} (see Fuchs \cite{Fu}) on slender groups, the
first `large' reflexive abelian groups are free groups or
(cartesian) products of $\Z$ - assuming for a moment that all
cardinals under consideration are $< \aleph_m$, the first
measurable cardinal. Also the members of the class of groups
generated by $\Z$ and taking direct sums and products
alternatively are reflexive, called Reid groups. Using a
generalized `Chase Lemma', which controls homomorphisms from
products of modules into direct sums of modules, Dugas and
Zimmermann-Huisgen \cite{DZ} showed that the class of Reid groups
is `really large'. Nevertheless there are more reflexive groups -
Eda and Otha \cite{EO} applied their `theory of continuous
functions on $0$-dimensional topological spaces' to find reflexive
groups not Reid-groups. As a by-product we also get dual groups
which are not reflexive, see also \cite{EM}. All these groups $G$
have the property that they are either free of finite rank or
\begin{eqnarray} \label{noniso}
\Z \oplus G \cong G.
\end{eqnarray}

As indicated in the abstract, we applied $\Diam_{\aleph_1}$ to
find examples $G$ of size $\aleph_1$, where (\ref{noniso}) is
violated. The obvious question to replace $\Diam_{\aleph_1}$ by
CH was the main goal of this paper. The question whether the Main
Theorem \ref{mainth} holds in any model of ZFC remains open. On
the other hand we are able to show that the conclusion of the
Main Theorem \ref{mainth} also follows in models of ZFC and
Martin's axiom MA, see \cite{GS1}. Hence  CH is not necessary to
derive the existence of these reflexive modules.

\section{Free modules with bilinear form and partial dual maps}

\begin{definition}\label{trip1}
Let $(\Phi, \FF_0, \FF_1)$ be a triple of a bilinear map $\Phi
:F_0 \oplus F_1 \arr R$ for some countable, free $R$-modules $F_i$
of infinite rank, $\FF_i \subseteq F_i^*  \  (i = 0,1)$ families
of dual maps subject to the following conditions
\begin{enumerate}
\item $\Phi$ is not degenerated. This is to say if
$\Phi(e,\ \ ) \in F_1^*$ or $\Phi(\ \ ,f) \in F^*_0$ is the
trivial map then $e = 0$ or $f = 0$, respectively. \item $\Phi$
preserves purity, that is $\Phi(e, \ ) \in_* F_1^*$ if $e \in_*
F_0$ and dually $\Phi(\  , f) \in_* F_0^*$ if $f \in_* F_1.$
\item $\FF_i$ is a countable, non-empty family of homomorphisms
$\varphi: F_\varphi \arr R$ such that $\Dom \varphi = F_\varphi
\subseteq_*F_i$. The set $\Dom \FF_i = \{F_\varphi : \varphi \in
\FF_i\}$ is well-ordered by inclusion, for $i \in \{0,1\}.$
\item For any $0 \neq x \in F_1$ and any finite subset $E \subset
\FF_0$ we have $\ker E \not\subseteq \ker \Phi(\  , x)$, and
dually for any $0 \neq y \in F_0$ and any finite subset $E \subset
\FF_1$, we have $\ker E \not\subseteq \ker \Phi(y, \  ).$
\end{enumerate}
\end{definition}

Here, and in the subsequent parts we use the following
\begin{notation} \label{trip2}
\begin{enumerate}

\item
$\FF$ is the collection of all triples $(\Phi, \FF_0, \FF_1)$ as
in Definition \ref{trip1}.
\item
If $E \leq \Hom (G,H)$ then $\ker E = \bigcap\limits_{\varphi \in
E} \ker \varphi.$
\item Similarly $\ker \Phi(\ , E) =
\bigcap\limits_{e \in E} \ker \Phi(\ ,e)$ and dually.
\end{enumerate}
\end{notation}

Next we define a partial order on $\FF$.
\begin{definition}\label{order}\begin{enumerate}
\item $(\Phi, \FF_0, \FF_1) \subseteq (\Phi', \FF_0', \FF'_1) \iff $
\item $\begin{cases}
\mbox{(a) }\ \Phi \subseteq \Phi' \mbox{ and } \Dom \Phi
\subseteq_* \Dom \Phi' \mbox{ is a pure submodule.}\\ \mbox{(b)
 } \mbox{ If } \varphi \in \FF_i \mbox{ then there is a unique }
\varphi' \in \FF'_i \mbox{ such that } \varphi \subseteq \varphi'
\mbox{ and }\\ \Dom \varphi \subseteq_* \Dom \varphi'.

\end{cases}$
\end{enumerate}
\end{definition}

We will construct the reflexive modules of size $\aleph_1$ by
using an order preserving continuous map and let
\begin{eqnarray} \label{map}
p: {}^{\omega_1 > }2 \arr \FF\ \ (\eta \arr p_\eta)
\end{eqnarray}
from the tree ${\bf T}$ of all branches of length $\lg(\eta) =
\alpha$ $$\eta : \alpha \arr 2 = \{0,1\}\  \mbox{for all}\ \alpha
< \omega_1.$$

 The order on $\bf T $ is defined naturally by
extensions, i.e. if $\eta, \eta' \in $ {\bf T}, then $\eta \leq
\eta'$ if and only if $\eta \subseteq \eta'$ as maps. Hence $\lg
(\eta)\leq\lg (\eta')$ and $\eta'\restr\lg (\eta) = \eta$, and we
will require that $p_\eta \subseteq p_{\eta'}$ by the ordering of
$\FF$ as defined in Definition \ref{order}. If $\eta \in
{}^{\omega_1} 2$, then $p_{\eta \restr \alpha} \ \ (\alpha \in
\omega_1)$ is linearly ordered and the triple
\begin{eqnarray} \label{triple}
p_\eta = \bigcup\limits_{\alpha \in \omega_1}\ p_{\eta \restr
\alpha}
\end{eqnarray}
is well-defined. In details we have $ \Phi_\eta : F_{0 \eta}
\oplus F_{1 \eta} \arr R$ is defined by continuity such that
$\Phi_\eta = \bigcup\limits_{\alpha \in \omega_1} \Phi_{\eta
\restr \alpha}$ and $F_{i \eta} = \bigcup\limits_{\alpha \in
\omega_1} F_{i \eta \restr \alpha}$. The bilinear forms $\Phi_\eta
: F_{0 \eta} \oplus F_{1 \eta} \arr R$ will be our candidates for
modules $G$ as in the Main Theorem \ref{mainth}. First we will
show that $\FF \neq \emptyset$ and the arguments will be refined
for Lemma \ref{sub}.

\begin{lemma}\label{set} \  The partially ordered set $\FF$ is non-empty.
\end{lemma}

{\bf Proof.} Choose $F_0 = F_1 = \bigoplus\limits_{n \in \omega}
e_n R$ and extend $$\Phi(e_i, e_j) = \delta_{ij} = \left \{
\begin{array}{@{}rrr@{}} 0 & if & i \neq j\\ [3pt] 1 & if & i = j
\end{array} \right.$$
linearly to get a bilinear map $\Phi : F_0 \oplus F_1 \arr R.$ The
map $\Phi$ satisfies Definition \ref{trip1} for $\FF_0 = \FF_1
=\emptyset$. Next we want find $\FF_0 = \{\varphi_0\}$ and $\FF_1
= \{\varphi_1\}$. If $\varphi_i \  (i \in 2)$ satisfy Definition
\ref{trip1} $(iv)$ and $\Dom \varphi_i = F_i$ then obviously
$$(\Phi, \FF_0, \FF_1) \in \FF.$$ We will work for $\varphi =
\varphi_0 : F_0 \arr R$ and enumerate $F_1 \setminus \{0\}  = \{
x_i :\ i \in \omega\}$. If $\varphi_i = \Phi(\ , x_i)$ then
$\varphi_i \neq 0$ by $(i)$ and for $(iv)$ we must show that
\begin{eqnarray} \label{ker}
\ker \varphi \not\subseteq \ker \varphi_i\  \mbox{for all}\  i \in
\omega
\end{eqnarray}

Write $F_0 = \bigoplus\limits_{i \in \omega} L_i$ such that each
$L_i = e_i R \oplus e'_i R$ is free of rank 2 and let $L'_i =
\bigoplus\limits_{j \neq i} L_j$ be a complement of $L_i$ . If
$L_i \varphi_i = 0$ then we use $\varphi_i \neq 0$ to find some
$0\neq y \in_* L'_i $ such that $y \varphi_i \neq 0.$ Choose a new
complement of $L'_i$ and rename it $L_i = (e_i + y) R \oplus e'_i
R$.

Hence $L_i \varphi_i \neq 0$ and there is a pure element $y_i \in
L_i$ with $y_i \varphi_i \neq 0$. We found an independent family
$\{y_i : i \in \omega\}$ with $F_0 = \bigoplus\limits_{i \in
\omega} y_i R \oplus C\ $ for some $ 0 \neq C \subseteq F_0$.
Choose $\varphi \in \Hom (F_0, R)$ such that $y_i \varphi = 0$ for
all $i \in \omega$ and $\varphi \restr C \neq 0$. Hence $$y_i \in
\ker \varphi \setminus \ker \varphi_i \ \mbox{for all} \ i \in
\omega$$ and (\ref{ker}) holds. Hence Definition \ref{trip1} holds
and $\FF_1$ can be chosen dually. $\hfill \square$

\bigskip

The crucial step in proving the next result is again verification
of Definition \ref{trip1} $(iv)$, this time for $\Phi'(\ , x)$.
The proof is similar to the last one, hence we can be less
explicit and just refine the old arguments.

\begin{lemma}\label{sub}\  Let $p = (\Phi, \FF_0, \FF_1) \in \FF$ be with
$\Dom \Phi = F_0 \oplus F_1$ and $\varphi \in L^*$ for some pure
submodule $L$ of finite rank in $F_0$. Then we find $p \subseteq
p' = (\Phi', \FF_0, \FF_1) \in \FF$ with $\Dom \Phi' = F_0 \oplus
F'_1, \  F'_1 = F_1 \oplus x R$ and $\varphi \subseteq \Phi'(\ ,
x).$
\end{lemma}
{\bf Proof.} First we want to extend $\varphi$ to $\varphi': F_0
\arr R$ such that
\begin{eqnarray} \label{finite}
\ker E \not\subseteq \ker \varphi' \mbox{ for all finite } E \leq
\FF_0.
\end{eqnarray}
Enumerate all $\ker E$ for $E \leq \FF_0$ finite by $\{K_i : \  0
\neq i \in \omega\}$.
 Hence $K_i = \bigcap\limits_{\varphi \in E_i} \ker\varphi$
is a pure submodule of $F_0$ with $F_0 / K_i$ free of rank $\leq
|E_i|$. As in the proof of (\ref{sub}) we can choose inductively
$\bigoplus\limits_{i \in \omega} L_i = F_0$ and $0 \neq y_i \in _*
L_i \cap K_i$  for $i > 0$ and $L_0 = L.$ For some $C \subseteq
F_0$ we have $F_0 = \bigoplus\limits_{i \in \omega} y_i R \oplus L
\oplus C$. Now we extend $\varphi \in L^*$ such that $y_i \varphi'
= 1$ and $\varphi \restr C = 0.$ Hence $$y_i \in K_i \setminus
\ker \varphi' \mbox{ for all }\ 0 \neq i \in \omega$$ and
(\ref{finite}) holds. Finally we extend $\Phi$ to $\Phi' : F_0
\oplus (F_1 \oplus x R) \arr R$ by taking $\Phi'(\ , x) =
\varphi'$. Condition (\ref{finite}) carries over to $$\ker E
\not\subseteq  \ker \Phi(\ , x) \mbox{ for all finite }  E \leq
\FF_0.$$ From this it is immediate that $ \ker E \not\subseteq
\ker \Phi(\ , y)$ for all $0 \neq y \in F'_1$  and finite  sets $
E \leq \FF_0$. Hence Definition \ref{trip1}$(iv)$ holds, $p \leq
(\Phi', \FF_0, \FF_1) \in \FF$ and $\varphi\subseteq\Phi(\ , x)$.
$\hfill\square$
\bigskip

Next we move elements from $\Hom (F_\varphi, R)$ for some $\varphi
\in \FF_i$ to $\FF_i$.

\begin{lemma} \label{PSI}\  Let $p = (\Phi, \FF_0, \FF_1) \in \FF$ be with $\Dom
\Phi = F_0 \oplus F_1$. If $\psi \in \Hom(F_\varphi, R)$ for some
$\varphi \in \FF_0$ such that $$\ker \psi \cap \ker E
\not\subseteq \ker \Phi(\ , x)\ \mbox{for all finite}\  E \leq
\FF_0\ \mbox{and}\  x \in F_1,$$ then $$p \leq p' = (\Phi, \FF'_0,
\FF_1) \in \FF \mbox{ for }  \FF'_0 = \FF_0 \cup \{\psi\}$$
\end{lemma}
{\bf Proof.}\  We only have to check Definition \ref{trip1}$(iv)$
for finite subsets of $\FF'_0 = \FF_0 \cup \{\psi\}.$ But this
follows by hypothesis on $\psi$.
\bigskip

The proof of the following observation is obvious.

\begin{observation} \label{chain}\
If $p_n = (\Phi_n, \FF_{0n}, \FF_{1n})\ (n \in \omega)$ is an
ascending chain of elements $p_n \in \FF$ and elements in $\FF_i$
are unions of extensions in $\FF_{in} \subseteq \FF_{in+1} \ (n
\in \omega)$, then $p = (\bigcup\limits_{n \in \omega} \Phi_n,
\FF_0, \FF_1) \in \FF.$
\end{observation}

\begin{definition} \  If $(\Phi, \FF_0, \FF_1) \in \FF$ with $\Dom \Phi =
F_0 \oplus F_1$ then $\varphi \in F_0^*$ is essential for $\Phi$
if for any finite rank summand $L$ of $F_0$ and any finite subset
$E$ of $F_1$ there is $g \in F_0 \setminus L$ with $g \varphi \neq
0$ and $\Phi(g, e) = 0$ for all  $e \in E$.
\end{definition}

The notion `$\varphi \in F_1^*$ is essential for $\Phi$ ' is dual.

If $g \varphi = 0 = \Phi(g, e)$ for $e \in E$ and some $E \leq
\FF_1$ then $g \in \langle \Phi(\ , e)  : e \in E \rangle
\subseteq F^*_0$ by induction on $|E|$. Hence $\varphi \in F_0^*$
is essential for $\Phi$ is equivalent to say that $\varphi$ is not
in $\Phi(\ ,F_1)$ modulo summands of finite rank. This leads to
the following

\begin{observation} \label{ines}  Let $(\Phi, \FF_0, \FF_1) \in \FF$ with
$\Dom \Phi = F_0 \oplus F_1$. If $E \leq F_1$ is a finite subset
and $\varphi \in F^*_0$ inessential for $\Phi$ with  $\ker \varphi
\subseteq \ker E$ i.e. $$(\Phi(x, e) = 0 \mbox{ for all }\ e \in
E) \Longrightarrow x \varphi = 0,$$ then there is $e_0 \in \langle
E \rangle \subseteq F_1$ such that $\varphi = \Phi(\ ,e_0)$.
\end{observation}

{\bf Proof.} \ By induction on $|E|$.

\begin{1killlemma} \label{Dom}
\  Suppose $\varphi \in F^*_0$ is essential for $\Phi$ with $\Dom
\Phi = F_0 \oplus F_1$ and $p = (\Phi, \FF_0, \FF_1) \in \FF$.
Then we find $p \leq p' \in \FF$ with $p' = \  (\Phi', \FF_0',
\FF_1'),\  \Dom\Phi' = F'_0 \oplus F_1$ and $F'_0 = \langle F_0, y
\rangle \subseteq_* \F_0$ for some $y \in \F_0$ such that
$\varphi$ does not extend to $\varphi' : F'_0 \arr R$.
\end{1killlemma}

{\bf Remark} A dual lemma holds for $\varphi \in F^*_1$.
\bigskip

{\bf Proof.}  Let $F_0 = \bigoplus\limits_{i \in \omega} e_i R$
and $F_1 = \bigoplus\limits_{n \in \omega} f_n R$. First we apply
that $\varphi$ is essential. It is easy to find inductively
elements $g_n \in F_0 \setminus F^n$ with $F^n = \langle e_i, g_i
: i < n \rangle_* \subseteq F_0$ such that the following holds
\begin{enumerate}
\item $\Phi(g_n, f_i) = 0$ for all $i < n.$
\item $\bigoplus\limits_{i< n} g_i R$ is a direct summand - also
$\bigoplus\limits_{i \in \omega}\  g_i R$ is a summand of $G$.
\item $g_n \varphi \neq 0$ for all $n \in \omega$.
\end{enumerate}

Decompose $\omega$ into a disjoint union of infinite subsets $S_i
\ (i \in \omega)$ and let $\{E_i: 0 \neq i \in \omega\}$ be an
enumeration of all finite subsets of $\FF_1$ and write $K_i = \ker
E_i$ for all $i > 0$. In order to get $K_i \not\subseteq \ker \Phi
(x,\ )$ for all   $i > 0, x \in F'_0 \setminus F$, we choose for
each $n \in S_i$ an element $k_n \in K_i$ such that $\Phi(g_n,
k_n) \neq 0$. This is possible as $F_1 / K_i$ is free of finite
rank $\leq |E_i|$, hence $K_i$ is `quite large'.

Now we use $g_n \varphi \neq 0\  (n \in \omega)$ and apply
Observation \ref{obs1} (for $S_0$) and choose a  sequence
$\epsilon_n \in \{0,1\} ( n \in \omega)$ and suitable $q_n\in S $
as in (\ref{qn-s}) such that
\begin{eqnarray}\label{r}
r = \sum\limits_{n \in \omega} (g_n \varphi)q_n  \epsilon_n \in \R
\setminus R
\end{eqnarray}
\begin{eqnarray}\label{not0}
\mbox{ If } n \in S_i, \mbox{ then } \sum^n\limits_{j = 0}  \Phi
(g_j, k_i) q_j \epsilon_j  \not\equiv 0 \ \mod\ q_{n+1}
\end{eqnarray}

If $s \in \omega$ and $\epsilon_s = 1$, then let $$y_s =
\sum\limits_{s \leq n \in \omega} (q_s)^{-1} q_n g_n \epsilon_n
\in \F_0$$ and consider the $R$-module $F'_0 = \langle F_0, y_0
\rangle_* \subset \F_0$, which can be generated by $$F'_0 =
\langle F_0, y_s R : \ s \in \omega, \epsilon_s = 1\rangle.$$

Note that $F'_0$ is a countable $R$-module. It is easy to see that
$F'_0$ is free - either apply Pontryagin's theorem (Fuchs \cite
[p. 93]{Fu}) or determine a free basis. The bilinear form $\Phi:
F_0 \oplus F_1 \arr R$ extends uniquely to $\Phi' : F_0' \oplus
F_1 \arr \widehat{R}$ by continuity and density. We want to show
that
\begin{eqnarray} \label{Phi}
(\Phi', \FF_0, \FF_1) \in \FF.
\end{eqnarray}

First we claim that \Im $\Phi' \subseteq R$. By (i) we have
$$\Phi'(y_0, f_j) = \Phi'(\sum\limits_{n \in \omega} \epsilon_n
q_n g_n, f_j) = \sum\limits_{n \in \omega} \Phi(g_n,
f_j)q_n\epsilon_n = \sum\limits_{n < j}\Phi(g_n, f_j)q_n\epsilon_n
\in R$$ and $\Im \Phi' \subseteq R$ follows.

We also must check $(iv)$ from Definition \ref{trip1} for the new
elements $y \in F'_0 \setminus F_0$. It is enough to consider
$$\ker E \not\subseteq \ker \Phi(y_s, \ ) \mbox{ for all } s \in
\omega \mbox{ with } \epsilon_s = 1.$$ By definitions and
enumerations this is equivalent to say that $$K_i \not\subseteq
\ker\Phi(y_s,\ )\ \mbox{for each}\  i > 0.$$ If $n \in S_i$ and $n
> s$ by (\ref{not0}) we have that $$\Phi(y_s, k_i) =
\sum\limits_{j \in \omega} \Phi(g_j, k_i) q_j \epsilon_j
\equiv\sum^n\limits_{j=0} \Phi(g_j, k_i) q_j \epsilon_j\not\equiv
0\ \mod q_{n+1}$$ Hence $k_i \not\in \ker\Phi(y_s, \ )$ but $k_i
\in K_i$ and (\ref{Phi}) follows. Finally we must show that
$\varphi \in F^*_0$ does not extend to $(F'_0)^*$. By continuity
$\varphi : F_0 \arr R$ extends uniquely to $\varphi' : F'_0 \arr
\R$. However $$y_0 \varphi = (\sum\limits_{n \in \omega}
g_nq_n\epsilon_n)\varphi' = \sum\limits_{n \in \omega} (g_n
\varphi)q_n\epsilon_n = r \in \R \setminus R $$ by (\ref{r}),
hence $\varphi$ does not extend to $F'_0 \arr R$. $\hfill\square$

\begin{2killlemma} \label{2kill}
\  Let $p = (\Phi, \FF_0, \FF_1) \in \FF$ and $\eta: F_\varphi
\arr F_\varphi$ be some monomorphism with $F_\varphi = x_0 R
\oplus F_\varphi\eta$ for some $\varphi \in \FF_0$ with $x_0 \in
F_\varphi = \ker \varphi$. Then there is $p \leq p' = (\Phi',
\FF'_0, \FF'_1) \in \FF$ with $\Dom \Phi' = F'_0 \oplus F'_1,
\varphi \subseteq \varphi' \in \FF'_0$ such that $\eta$ does not
extend to a monomorphism $$\eta'' : F_{\varphi'} \arr F_{\varphi'}
\subset_* F''_0
 \mbox{ with } F'' _{\varphi'} = x_0 R \oplus F''_{\varphi'} \eta''$$

where $p' \leq p'' := (\Phi'', \FF_0'', \FF_1'')$ and $\Dom \Phi''
= F{''}_0 \oplus F''_1$.
\end{2killlemma}
{\bf Proof.}\ In order to satisfy Definition \ref{trip1} $(iv)$
for the new $p' \in \FF$ we argue similar to the First Killing
Lemma \ref{Dom}. Let $\omega = \bigcup\limits_{i \in \omega} S_i$
be a decomposition into infinite subsets $S_i$ and $\{K_i : 0 \neq
i \in \omega\}$ be all kernels $K_i = \ker E_i$ for an enumeration
of finite subsets $E_i$ of $\FF_1$. The set $S_0$ will be used for
killing $\eta$ and the $S_i \  (i > 0)$ are in charge of $K_i$ and
$(iv)$ above. Extending $\Phi \subset \Phi'$ we must ensure $\Im
\Phi' \subseteq R$. Hence we construct an increasing sequence $s_n
\in \omega \ (n \in \omega)$ and pose more conditions on $s_n$
later on. If $F' = \Dom \varphi$ then $F' = x_0 R \oplus F'\eta$
and if $x_i = x_0 \eta^i$ we get
\begin{eqnarray}\label{new}
F' = \bigoplus\limits_{i<n} x_i R \oplus F' \eta^n  \mbox{ for all
}  n \in \omega.
\end{eqnarray}

Note that $\bigoplus\limits_{i<n} x_i R$ is pure of finite rank,
hence $F_0 = \bigoplus\limits_{i<n} x_i R\oplus C_n $ for some
$C_n \subseteq F_0$.  Now we construct $T_n =
\bigoplus\limits_{s_n \leq i < s_{n+1}} e_i R$ from $F_0 =
\bigoplus\limits_{i \in \omega} e_i R$ and refine an argument from
G\"obel, Shelah \cite{GS}. Obviously $F_0 = \bigoplus\limits_{n
\in \omega} T_n$ and let $e^n \in \Hom(F_0, R) \ (n \in \omega)$
be defined by
$$e_i e^n = \delta_{i,n} = \left \{
\begin{array}{@{}rrr@{}} 0 & if & i \neq n\\ [3pt] 1 & if & i = n
\end{array} \right.$$
for all $i \in \omega$. If $i < n$ then $\pi_i : F_0 \arr R$
denotes the projection modulo $\bigoplus\limits_{i \neq j<n} x_j R
\oplus C_n$, moreover let $F_1 = \bigoplus\limits_{j \in \omega}
f_j R$ be as before. We now seek for elements $w_n \in F_0 \  (n
\in \omega)$ subject to the following four conditions
\begin{enumerate}
\item $0 \neq w_n \in T_n$ and $w_n \eta \in T_n$.
\item $\Phi(w_n, f_n) = w_n e^k = w_n \eta e^k = 0$ for all $k < s_n$
\item If $n \in S_0,$ then let $\pi^*_n$ be the projection $\pi_i$ with
$i$ maximal such that $$w_n \pi_i = 0 \neq w_n \eta \pi_i, \ s_n
\leq i < s_{n+1}$$
\item If $n \in S_i, i > 0$, there is $y_n \in K_i$ such that $\Phi(w_n, y_n)
\ne 0.$
\end{enumerate}

Suppose $s_0, \dots, s_n, w_0 \dots, w_{n-1}$ are constructed and
we want to choose $w_n, s_{n+1}$. Then pick $s_n < s_{n+1}$ such
that $$\{x_{s_n}, \dots, x_{4 s_n +1}\} \subseteq
\bigoplus\limits_{i \leq s_{n+1}} e_i R.$$

We want to choose $w_n = \sum^{4 s_n}\limits_{i = s_n} x_i a^n_i$
for some $a^n_i \in R$.  If $w_n e^k = 0$ for $k < s_n$, then $w_n
\in T_n$. Moreover $$w_n \eta = \sum^{4 s_n}\limits_{i=s_n} x_i
a^n_i \eta = \sum^{4 s_{n}}\limits_{i = s_{n}} x_{i+1} a^n_i =
\sum^{4 s_n +1}\limits_{i = s_n +1} x_i a^n_{i-1}$$ and $x_n \eta
\in  T$ because $x_{4s_n +1} \in \bigoplus\limits_{i \leq s_
{n+1}} e_i R$ as well.

Hence $(i)$ follows provided $w_n \neq 0$ is generated by those
$x_i's$. The conditions $(ii)$ can be viewed as a system of $3
s_n$ homogeneous linear equations in $4 s_n + 1 - s_n = 3 s_n + 1$
unknowns $a^n_i \in R$. We find a non-trivial solution $w_n \neq
0$ by linear algebra. Hence $(i)$ and $(ii)$ hold. Condition
$(iii)$ follows by hypothesis on $\FF_1$ for $E_i$ and $K_i = \ker
E_i$. \ Condition $(iii)$ finally follows by the action of $\eta$
on $w_n \neq 0$ and the maximality of $i$ with $s_n < i \leq 4
s_n$ and $w_n \pi_i = 0$ for $\pi^*_n = \pi_i$. Hence $(i),\dots,
(iv)$ follow.

As in the proof of the First Killing Lemma \ref{Dom} inductively
we choose a strictly increasing sequence $m_j \in \omega \ (j \in
\omega)$. If $m_j$ is defined up to $j \leq n$ we must choose
$m_{n+1}$ large enough such that $$\sum\limits_{j \leq n}
q_{m_{j}} \Phi(w_j, y_n) \not\equiv 0 \mod q_{m_{n+1}}.$$

This needs inductively the hypothesis that $\sum\limits_{j \leq n}
q_{m_{j}} \Phi(w_j, y_n) \neq 0$. If $n+1 \in S_i$, then $\Phi
(w_{n+1}, y_{n+1}) \neq 0$ by $(iv)$ and we may assume
$$\sum\limits_{j \leq n+1} p^m_{j} \ \Phi(w_j, y_{n+1}) \neq 0$$
hence the inductive hypothesis follows and we can proceed. By
Observation \ref{obs1} and $(iii)$ we also find $\epsilon_j \in
\{0,1\} \ (j \in S_0)$ such that
\begin{eqnarray}\label{Pi}
\sum\limits_{j \in S_0} (w_j \eta) \pi^*_j q_{m_j}\epsilon_j \in
\R \setminus R.
\end{eqnarray}
\ \\ Now we are ready to extend $\Phi$. Choose new elements $$z_k
= \sum\limits_{j \geq k} w_j q_{m_{j}} (q_{m_{k}})^{-1} \in \F_0$$

for all $k \in \omega$. Hence the submodule $F'_0 = \langle F_0,
z_0 \rangle_* \subseteq \widehat{F}_0$ purely generated by adding
$z = z_0$ is generated by $$F'_0 = \langle F_0, z_k : k \in \omega
\rangle.$$

Again we see that $F'_0$ is a countable, free $R$-module. The map
$\Phi: F_0 \oplus F_1 \arr R$ by continuity extends uniquely to
$$\Phi' : F'_0 \oplus F_1 \arr \widehat{R}.$$

Recall that $F'_0 / F_0$ is S-divisible, hence $F_0$ is S-dense in
$F'_0$ in the S-adic topology. First we must show that $\Im \Phi'
\subseteq R$. We apply $(ii)$ and continuity to see that
$$\Phi'(z, f_k) = \Phi'(\sum\limits_{n \in \omega} q_{m_{n}} w_n,
f_k) = \sum\limits_{n \in \omega} q_{m_{n}} \Phi(w_n, f_k) =
\sum\limits_{n \leq k} q_{m_{n}} \Phi(w_n, f_k) \in R,$$ hence
$\Phi' : F'_0 \oplus F_1 \arr R$.

We also must show Definition \ref{trip1} $(iv)$ for the new
elements $z_t \in F'_0$. We have $K_i = \ker E_i$ and $S_i$ is
unbounded. Hence we find $n \in S_i, n > t$ and  $y_n \in K_i$
such that $$\Phi(w_n, y_n) \neq 0.$$ We apply $\Phi'$ to $(z_t,
y_n)$ and get $$\Phi'(z_t, y_n) \equiv \sum\limits_{t \leq j \leq
n} q_{m_j} (q_{m_t}^{-1} \Phi(w_j, y_n) \not\equiv 0 \mod
q_{m_{n+1}} q^{-1}_{n_t}$$

hence $\Phi(z_t, y_n) \neq 0$ and $y_n \in K_i$. This is
equivalent to say that $\ker E_i \not\subseteq \ker \Phi(z_t ,\ )$
and Definition \ref{trip1} $(iv)$ follows. Hence $(\Phi', \FF_0,
\FF_1) \in \FF$ with $\Dom\Phi' = F'_0 \oplus F_1.$

Next we extend $\Phi'$ under the name $\Phi'$. Let $F'_1 = F_1
\oplus f R$ be a free rank-1 extension. We want $\Phi' : F'_0
\oplus F'_1 \arr R$ and must define $\Phi'(\ , f) : F'_0 \arr R$.
Put $\Phi'(\ , f) \restr T_n = \epsilon_n \pi^*_n \restr T_n$ if
$n \in S_0$ and $\Phi'(\ , f) \restr T_n = 0$ otherwise. By linear
extension and $F_0 = \bigoplus\limits_ {n \in \omega} T_n$ the map
$\Phi'(\ ,f) : F_0 \arr R$ is well-defined. It extends further by
continuity to $$\Phi' : F'_0 \oplus F'_1 \arr \widehat{R}.$$ Again
we must show that $\Im\Phi'\subseteq R$. Note that $\Phi'(w_n, f)
= w_n \epsilon_n \pi^*_n = 0$ for $n \in S_0$ from (iii), and
$\Phi'(w_n, f) = 0$ for $n \in \omega \setminus S_0$ by the above,
hence $$\Phi'(z, f) = \sum\limits_{n \in \omega} q_{m_{n}}
\Phi'(w_n, f) = 0$$ and $\Im \Phi'\subseteq R$ follows.

We also must check condition $(iv)$ of Definition \ref{trip1} for
the new element $f \in F'_1$. Recall $K_i\ (i \in \omega)$ is a
list of all kernels $\ker E_i$ for finite sets in $\FF_0$. Also
enumerate all $T_n$'s with $n \in S_0$ and $r_n = 1$ as, say $T^i
= T_{n_{i}} \ (i \in \omega)$. Choosing the ranks of the $T^{i'}s$
large enough we can find (as before) $y'_i \in K_i \cap T^i$ with
$y'_i \pi^*_{n_{i}} \neq 0$. Then $$y'_i \in \ker E_i  \mbox{ but
} \Phi'(y'_i, f) = y'_i \pi^*_{n_{i}} \neq 0,$$  hence $y'_i
\not\in \ker\Phi(\ ,f)$ as desired for $(iv)$ above.

Finally we must get rid of $\eta$ by showing that there is no
extension $\eta'' \supset \eta$ as stated in the Lemma. Otherwise
we have $\Phi'\subset\Phi'' : F''_0 \oplus R\arr R$ and $$x_0R
\oplus F'' \eta'' = F''$$ for $\Dom\eta'' = F'' \subset_* F''_0$
with $z \in \Dom\eta'\subseteq F''$.
 Hence $r = \Phi'' (z \eta'', f) \in R$. On the other
hand $$r = \Phi'' (z \eta'', f) = \Phi'' (\sum\limits_{j \in
\omega} q_{m_{j}} (w_j \eta'), f) = \sum\limits_{j \in \omega}
q_{m_{j}} \Phi'(w_j \eta', f) = \sum\limits_{j \in S_{0}}
q_{m_{j}} \epsilon_j (w_j \eta \pi^*_j) \in \widehat{R} \setminus
R$$ is a contradiction. The lemma follows. $\hfill{\square}$

\section{Construction of reflexive modules assuming CH}
Let $\{S_0, \dots ,S_5\}$ be a decomposition of the set of all
limit ordinals in $\omega_1$ into stationary sets. Using CH we can
enumerate sets of cardinality $2^{\aleph_{0}}$ by any of these
stationary sets of size $\aleph_1$.\\ Let
\begin{eqnarray}\label{CH}
{\bf T} = {}^{\omega_{1 >}} 2 = \{\eta_\alpha : \alpha \in
\omega_1\}
\end{eqnarray}

be an enumeration of the tree with countable branches such that
\begin{eqnarray}\label{tree}
\eta_\alpha \subseteq \eta_\beta  \Longrightarrow  \alpha < \beta
\end{eqnarray}

We are working in the `universe' $\omega_1$ and let
\begin{eqnarray}
\Map (\omega_1, R) = \{\varphi_i : i \in S_0 \} = \{\varphi_i : i
\in S_1\} = \{\varphi_i : i \in S_2\} = \{\varphi_i : i \in S_3\}
\end{eqnarray}
be four lists of all partial functions $\varphi : \delta \arr R$
for $\delta \in \omega_1$ any limit ordinal with $\aleph_1$
repetitions $\varphi_i$ for each $\varphi$. Similarly let $$\Map
(\omega_1, \omega_1) = \{\mu_i : i \in S_4\} = \{\mu_i : i \in
S_5\}$$ be two lists of all partial maps $\mu: X \arr X \subseteq
\omega_1$ with countable domains $X$ given by the prediction Lemma
\ref{predict} for $E \in \{S_4, S_5\}$.

Inductively we want to construct an order preserving continuous
map $$p : {}^{\omega_{1}
>} 2 = {\bf T} \arr \FF,\ ( \eta \arr p_\eta)$$

subject to certain conditions $(a), \dots , (d)$ dictated by the
proof of the main theorem and stated below. Some preliminary words
are in order. Using the enumeration (\ref{CH}) we may write
\begin{eqnarray}\label{oplus}
p_{\eta_{\alpha}} = p_\alpha = (\Phi_\alpha, \FF_{0 \alpha},
\FF_{1\alpha}) \mbox{ and } \Dom\Phi_\alpha = F_{0 \alpha} \oplus
F_{1 \alpha}
\end{eqnarray}
If $p$ is order preserving then $(\eta_\alpha \subseteq \eta_\beta
\Longrightarrow p_\alpha \subseteq p_\beta)$  and continuity
applies. If $\eta \in {}^\alpha 2$ has length $\lg (\eta) =
\alpha$ and $\alpha < \omega_1$ is a limit ordinal, then
$\bigcup\limits_{\beta<\alpha} \eta \restr \beta = \eta$ and $\eta
\restr \beta \subseteq \eta$, hence
\begin{eqnarray}\label{eta}
\eta \in {}^\alpha 2 \Longrightarrow\  p_\eta =
\bigcup\limits_{\beta < \alpha} p_{\eta \restr \beta}
\end{eqnarray}

If all $p_{\eta \restr \beta} \in \FF$ then $p_\eta \in \FF$ by
Observation \ref{chain}. Hence we have no problem in defining $p$
at limit stages by continuity. It remains to consider inductive
steps for $p$. If $\eta, \eta' \in {}^{\omega_{1}>}2$ then $\gamma
= \br (\eta, \eta')$ denotes the branching point of $\eta, \eta'$,
this is to say that
\begin{enumerate}
\item $\gamma < \min \{\lg (\eta), \lg (\eta')\}$
\item $\eta \restr \gamma = \eta' \restr \gamma$
\item $\eta (\gamma) \neq \eta' (\gamma).$
\end{enumerate}
If $\eta$ is a branch of length $\alpha$ and $i \in \{0,1\}$ then
$\eta' = \eta\wedge \{i\}$ is a branch of length $\alpha + 1$ with
$\eta' \restr \alpha = \eta$ and $\eta' (\alpha) = i.$ Now we
continue defining $p : {}^{\omega_1 >} 2 \arr \FF$. Generally we
require
\begin{description}{\it
\item (a) \ If $\varphi \in \FF_{i \alpha}$ then $\Dom \varphi =
F_{i \beta}$ for some $\eta_\beta \subseteq \eta_\alpha.$
\item (b) The set $\{\Dom \varphi : \varphi \in \FF_{i \alpha}\}$
is a well ordered set of pure submodule of $F_{i \alpha}$.
\item
(c) \ If $\gamma = \br(\eta_\alpha, \eta_\beta)$ and $\eta_\alpha
\restr \gamma  = \eta_\beta \restr \gamma  = \eta_\epsilon$ for
some $\epsilon < \omega_1$, and $y \in F_{1 \beta} \setminus F_{1
\epsilon}$ then $\Phi^\beta (y,\ ) \restr F_{0 \epsilon} \in
\FF_{0 \alpha}.$ Dually, if $y \in F_{0 \beta} \setminus F_{1
\beta}$ then $\Phi^\beta (y, \ ) \restr F_{1 \epsilon} \in \FF_{1
\alpha}.$
\item  (d) \
If $\nu \in {}^\gamma 2$ is a branch of length $\gamma,\ i
\in\{0,1\} $ and $\eta = \nu\wedge \{i\}$ then we want to define
$p_\eta$ depending on $\gamma$.
\begin{enumerate}
\item For $\nu =\emptyset$ choose $p_{\{i\}}$ by Lemma \ref{sub}.
\item If $\gamma \in S_0$ then we want to enlarge $\Phi_\nu$
to make the evaluation map injective: If $\varphi_\gamma \restr
F_\varphi$ for some $\varphi \in \FF_{0 \nu}$ is a partial
$R$-homomorphism $F_\varphi \arr R$ with $\Dom(\varphi_\gamma
\restr F_\varphi)$ a pure $R$-submodule of finite rank of
$F_\varphi$, then we apply Lemma \ref{PSI} such that $p_\nu
\subset p_\eta$ with $F_{1 \eta} = F_{1 \nu} \oplus x_\eta R, \
F_{0 \eta} = F_{0 \nu}$, $\Phi_\nu \subset \Phi_\eta$ and
$\varphi_\gamma \subset \Phi_\eta (\ , x_\eta).$ If
$\varphi_\gamma$ does not satisfy the requirements, then we choose
any  $\Phi_\nu \subset \Phi_\eta$.
\item If $\gamma \in S_1$ then we argue as in $(ii)$ but dually. A dual
version of Lemma \ref{PSI} provides $\Phi_\nu \subset \Phi_\eta$
and $\varphi_\gamma \subset \Phi_\eta (x_\eta,\ )$ if
$\varphi_\gamma$ meets the requirements.
\item If $\gamma \in S_2$ then we want to kill bad dual maps to make
the evaluation map surjective. If $\varphi_\gamma \restr F_{0
\nu}$ is an $R$-homomorphism $F_{0 \nu} \arr R$ which is essential
for $\Phi_\nu$ then we apply the First Killing Lemma \ref{Dom} to
find $p_\gamma \subset p_\eta$ such that $F_{0 \eta} = \langle
F_{0 \nu}, y \rangle \subseteq_* \widehat{F}_{\nu 0}$ for some $y
\in \widehat{F}_{\nu 0}$ and $\varphi_\gamma \restr F_{0 \nu}$
does not extend to $F_{0 \eta} \arr R.$
\item If $\gamma \in S_3$ and $\varphi_\gamma \restr F_{1 \nu}$ is an
essential $R$-homomorphism for $\Phi_\nu$ then we argue as in
$(iv)$ but dually.
\item If $\gamma \in S_4$ then we want to get rid of potential
monomorphisms $\eta$ of the final module $G$ with $G \eta \oplus
xR = G.$ If $\varphi \in \FF_{0\nu}, \ \mu = \varphi_\gamma \restr
F_\varphi $ and $\mu: F_\varphi \arr F_\varphi$ is an
$R$-monomorphism such that $F_\varphi = x_\varphi R \oplus
F_\varphi \mu$, then we apply the Second Killing Lemma \ref{2kill}
to find $p_\nu \leq p_\eta $ such that $\mu$ does not extend to a
monomorphism $\mu'$ of an extension $F'_0$ of $F_{0\eta}$ with
$$p_\eta \subseteq p' = (\Phi', \FF'_0, \FF'_1),\ \Dom \Phi' =
F'_0 \oplus F'_1  \mbox{ and } F'_0 = x_\varphi R \oplus F'_0
\mu'.$$
\item If $\gamma \in S_5$, then we argue dually for some partial monomorphism $\mu$
with domain and range some $F_\varphi \subseteq_* F_{1\nu}$.
\item If $\gamma \in \omega_1$ is not a limit ordinal, then we are free to
choose any `trivial' extension $p_\nu \subset p_\eta$.
\end{enumerate}}
\end{description}

\section{Proof of the Main Theorem}
Recall from Section 3 that we are given an order preserving
continuous map $$p: \mbox{ \bf T } \arr \FF.$$ By continuity we
may extend $p$ to $$p_\eta = (\Phi_\eta, \FF_{0 \eta}, \FF_{1
\eta}) \mbox{ for all } \eta \in {}^{\omega_1}2.$$

It is immediate that
\begin{eqnarray}\label{a-free}
F_{i \eta} \mbox{ is } \aleph_1\mbox{-free  and } \Phi_\eta : F_{0
\eta} \oplus F_{1 \eta} \arr R \mbox{ is a bilinear form, }
\end{eqnarray}
where  $i \in \{0,1\},\ \eta \in {}^{\omega_1}2.$ Moreover,
$\Phi_\eta$ preserves purity and is not degenerated in the sense
of Definition \ref{trip1}. $\Phi_\eta : F_{0 \eta} \oplus F_{1
\eta} \arr R$ is our candidate for a reflexive modules, expressed
in an unusual way.

To see that  $\Phi_\eta$ preserves purity we must show that each
pure element $e \in_* F_{0 \eta}$ induces $\Phi_\eta (e,\ ) \in_*
F^*_{1 \eta}$ (and dually). We may restrict to the first case. If
$e \in_* F_{0 \eta}$ then $e \in_* F_{0 \eta \restr \alpha}$ for
any $\alpha \in \omega_1$ large enough, and $\Phi_{\eta \restr
\alpha} (e, \ ) \in_* F^*_{1 \eta \restr \alpha}$ by $(\Phi_{\eta
\restr \alpha}, \FF_{0 \eta \restr \alpha}, \FF_{1 \eta \restr
\alpha}) \in \FF$ and Definition \ref{trip1} $(ii)$. It follows
that $\bigcup\limits_ {\alpha \in \omega_{1}} \Phi_{\eta \restr
\alpha} (e,\ ) \in_* F^*_{1 \eta}.$

To see that $\Phi_\eta$ is not degenerated we consider $0 \neq e
\in F_{0 \eta}$. Hence $e \in e'R \subseteq_* F_{0 \eta \restr
\alpha}$ for a pure element $e'$ and any $\alpha \in \omega_1$
large enough. The partial homomorphism $\varphi: e'R \arr R$
defined by $e'\varphi = 1$ has a number $i \in S_0$ in the list
 and $\varphi = \varphi_i$ . By construction there is
 $y \in F_{1 \eta \restr \gamma}$ with
 $$\varphi \subset \Phi_{\eta \restr \gamma}
(\ , y).$$ Hence $0 \neq e \varphi = \Phi_{\eta \restr \gamma}(e,
y) = \Phi_\eta (e, y)$ and $\Phi_\eta$ is not degenerated.
$\hfill{\square}$

\begin{definition}\label{fullrep}
We will say that $(\Phi, F_0, F_1)$ with $\Dom \Phi = F_0 \oplus
F_1$ is fully represented if $\Phi_\eta (F_{0 \eta}, \ ) = F^*_{1
\eta}$ and $\Phi_\eta (\ , F_{1 \eta}) = F^*_{0 \eta}$.
\end{definition}

We claim that
\begin{eqnarray}\label{f-rep}
(\Phi_\eta,F_{0 \eta},F_{1 \eta}) \mbox{ is fully represented for
almost all } \eta \in {}^{\omega_1}2
\end{eqnarray}
There are at most $\eta \in W \subseteq  {}^{\omega_1}2$
exceptions with $|W| < 2^{\aleph_1}$.

Suppose for contradiction that $|W| = 2^{\aleph_1}$ and
\begin{eqnarray}\label{inW}
\varphi_\eta \in F^*_{1 \eta} \setminus \Phi_\eta (F_{0 \eta}, \
)\ \mbox {for all}\ \eta \in W
\end{eqnarray}
By a pigeon hole argument there are $\eta, \eta' \in W$ with
$\br(\eta, \eta') = \alpha$ and

(a) $\varphi_\eta \restr F_{1 \alpha} = \varphi_{\eta'} \restr
F_{1 \alpha}$

(b) $\varphi_\eta : F_{1 \eta} \arr R$, and $\varphi_{\eta'} :
F_{1 \eta'} \arr R$ are not represented.

Let $\psi = \varphi_\eta \restr F_{1 \alpha} = \varphi_{\eta'}
\restr F_{1 \alpha} : F_{1 \alpha} \arr R$ and recall that there
is some $\gamma \in S_3$ with  $\psi = \varphi_\gamma$ and
$(d)(v)$ of the construction applies. Hence $\psi$ is inessential
for $\Phi_\gamma$. There is a finite set $E \leq F_{0 \gamma}$
such that $$\Phi_\gamma (e, x) = 0\ \mbox{for all}\ e \in E,\ x
\in \Dom \psi = F_{1 \alpha} \mbox{ implies } x \psi = 0.$$ By
Observation \ref{ines} we have some $e \in \langle E \rangle
\subseteq F_{0 \gamma}$ such that $\psi = \Phi_\gamma(\ ,e)\restr
F_{1\alpha}$. The same argument applies for $\varphi_{\eta'}$ and
there are $\gamma'\in S_3$ and $e' \in F_{0 \gamma'}$ such that
$\gamma < \gamma'$ and $\psi = \Phi_{\gamma'}(\ ,e')\restr
F_{1\alpha}$. Hence $\psi = \Phi_{\gamma'}(\ , e') \restr F_\alpha
= \Phi_\gamma(\ , e) \restr F_\alpha$.

Finally we apply $(c)$ of the construction to get $\psi \in \FF_{0
\gamma'}$. Now Definition \ref{trip1} $(iv)$ applies and $\ker\psi
\not\subseteq \ker \Phi_{\gamma'}(\ , e')$ is a contradiction
because $\psi = \Phi_{\gamma'}(\ , e')$. The claim (\ref{f-rep})
follows. $\hfill{\square}$
\bigskip

Note that we did not use $(vi)$ so far. Hence without the Second
Killing Lemma \ref{2kill} we are able to derive reflexivity of the
modules $G_\eta$, which we will do next.

We will use the following notations.

$$\mbox{Let } I' = \{ \eta \in {}^{\omega_1} 2 \mbox{ such that }
p_\eta \mbox{ is not fully represented\} and } I = {}^{\omega_1}
2\setminus I'$$

\relax From (\ref{f-rep}) we see that $|I'| < 2^{\aleph_1}$, hence $|I| =
2^{\aleph_1}$ . If $\eta \in I$ and $i \in \{0,1\}$ then we also
fix the evaluation map $$\sigma_{i \eta} : F_{i \eta} \arr
F_{i\eta}^{**}$$ and claim that
\begin{eqnarray}\label{injectiv}
\sigma_{i \eta} : F_{i \eta} \arr F_{i \eta}^{**}\ \mbox{ is
injective for all } \eta \in I, \ i \in \{0,1\}.
\end{eqnarray}
{\bf Proof. }  We consider $\sigma= \sigma_{0 \eta}$ and apply
that $\Phi_\eta$ is not degenerated. If $0 \neq x \in F_{0 \eta}$
there is $y \in F_{1 \eta}$ such that $\Phi_\eta(x,y) \neq 0$.
Hence $\varphi := \Phi_\eta (\ , y) \in F_{0 \eta}^*$ and $x
\varphi = \Phi_\eta (x,y) \neq 0$, thus $x \sigma \neq 0$ and
$\sigma$ is injective. The case $\sigma_{1\eta}$ is similar.
$\hfill{\square}$

Next we show that
\begin{eqnarray}\label{surject}
\sigma_{i \eta} : F_{i \eta} \arr F^{**}_{i \eta} \mbox{ is
surjective for all }  \eta \in I, \  i \in \{0,1\}.
\end{eqnarray}

{\bf Proof.} First note that
\begin{eqnarray}\label{biject}
\Phi^\bullet_\eta : F_{0 \eta} \arr F^*_{1 \eta}\ \ (x \arr
\Phi_\eta(x,\ )) \mbox{ is bijective }
\end{eqnarray}
and
\begin{eqnarray}\label{2biject}
{}^\bullet \Phi_\eta : F_{1 \eta} \arr F^*_{0 \eta}(y \arr
\Phi_\eta\ (\ , y)) \ \ \mbox{ is bijective }
\end{eqnarray}
because  $\Phi_\eta$ is not degenerated and $\eta \in I$. Hence we
can identify $F_{1 \eta}$ and $F^*_{0\eta}$ by $^\bullet
\Phi_\eta$ and  $F^*_{0 \eta} = \Im (^\bullet \Phi_\eta) =
\Phi_\eta(\ , F_{1 \eta})$. Moreover $$F^{**}_{0 \eta} = (F^*_{0
\eta})^* = (F_{1 \eta})^* = \Im \Phi^\bullet_\eta = \Phi_\eta\
(F_{0 \eta},\ )$$ and for any $\varphi \in F^{**}_{0 \eta}$ we
find $f \in F_{0 \eta}$ with $\varphi = \Phi_\eta(f,\ )$.  We
consider the case $\sigma = \sigma_{0 \eta}$, and get
$$\Phi_\eta(\ , x)\sigma(f) = \Phi_\eta(f,x) = \Phi_\eta(\ ,x)
\Phi_\eta(f,\ )$$ for all $x \in F_{1 \eta}$ and $\Phi_\eta(\ ,x)$
runs through all of $F^*_{0 \eta}$. We derive $\sigma(f) =
\Phi_\eta(f, \ ) = \varphi$ and $\sigma$ is surjective. The case
$\sigma_{1 \eta}$ is similar. $\hfill{\square}$

\bigskip

We have an immediate corollary from (\ref{injectiv}) and
(\ref{surject}).
\begin{corollary} If $\eta \in I$ and $i \in \{0,1\}$ then
 $F_{i \eta}$ is a reflexive $R$-module. Moreover $I$ is a subset
 of $\ {}^{\omega_1} 2$ of cardinality $2^{\aleph_1}$.
\end{corollary}
For the {\bf Proof of the Main Theorem \ref{mainth}} we finally
must show that
\begin{eqnarray}\label{incong}
F_{i\eta} \not\cong R \oplus F_{i\eta} \mbox{ for any } \eta \in
I,\ i \in \{0,1\}.
\end{eqnarray}
We consider $F_{i\eta}$ with some monomorphism $\xi: F_{i\eta}
\arr F_{i\eta}$ such that $F_{i\eta} = F_{i\eta}\xi \oplus xR$ for
any $\eta \in I, i \in \{0,1\}$. By a back and forth argument
there is an $\alpha < \omega_1$ such that $F_{i\eta\restr\alpha} =
F_{i\eta\restr\alpha}\xi \oplus xR$ and we take $\psi = \xi\restr
F_{i\eta\restr\alpha}$ into consideration. There is some $\gamma
\in S_4$ such that $\psi = \varphi_\gamma$ and $\varphi_\gamma$ is
discarded by the construction. Hence $\xi$ does not exist.
$\hfill{\square}$
\bigskip

We would like to add a modification of our main result which can
be shown using one more stationary subset $S_6$ after introducing
inessential endomorphisms for our category of reflexive modules.
Let $\Fin(G)$ be the ideal of all endomorphisms $$\{\sigma \in
\End(G) : \ G\sigma \mbox{ has finite rank}\}$$ for some
torsion-free $R$-module $G$. Then we can find a `Killing-Lemma' in
Dugas, G\"obel \cite{DG}, see also \cite{CG}, which `takes care'
of all endomorphisms which are not in $\Fin(F_{i\eta})$ with $i,
\eta$ as above. Hence we can strengthen our Main Theorem
\ref{mainth} and get with slight modification from known results
the following

\begin{corollary}\label{ring}
(ZFC + CH) \  Let $R$ is a countable domain  but
not a field and $A$ be a countable $R$-algebra with free additive
structure $A_R$. Then there is a family of $2^{\aleph_1}$
pair-wise non-isomorphic reflexive $R$-modules $G$ of cardinality
$\aleph_1$ such that $G \not\cong R \oplus G$ and $\End(G) =
A\oplus \Fin(G)$ a split extension.
\end{corollary}

In particular $\End(G)/\Fin(G)\cong A$ and if $A$ has only trivial
idempotents like $A=R$ then $G$ in the Corollary \ref{ring} is
reflexive, essentially indecomposable of size $\aleph_1$ and does
not decompose into $G \cong R \oplus G$.

\noindent
R\"udiger G\"obel \\ Fachbereich 6, Mathematik und
Informatik
\\ Universit\"at Essen, 45117 Essen, Germany \\ {\small e--mail:
R.Goebel@Uni-Essen.De}\\ and \\ Saharon Shelah \\ Department of
Mathematics\\ Hebrew University, Jerusalem, Israel
\\ and Rutgers University, Newbrunswick, NJ, U.S.A \\ {\small
e-mail: Shelah@math.huji.ae.il}


\begin{thebibliography}{99}\label{lit}
\markright{}

\bibitem{B}  H. Bass, Finitistic dimension and a homological
generalization of semi-primary rings, Transact. Amer. Math. Soc.
{\bf 95}, 466 -- 488 (1960).

\bibitem{CG}  A. L. S. Corner, R. G\"obel, Prescribing
endomorphism algebras - a unified treatment, Proc. London Math.
Soc. (3) {\bf 50}, 447 -- 479 (1985).

\bibitem{DS}  K. Devlin, S. Shelah, A weak version of $\Diam$
which follows from $2^{\aleph_0} < 2^{\aleph_1}$, Israel J. Math.
{\bf 6}, 239 -- 247 (1978).

\bibitem{DG} M. Dugas, R. G\"obel, Endomorphism rings of separable
torsion-free abelian groups, Houston J. Math. {\ 11}, 471 -- 483.

\bibitem{DZ}  M. Dugas, B. Zimmerman-Huisgen, Iterated direct sums and
products of modules, pp. 179 -- 193 in `Abelian Group Theory',
Springer LNM {\bf 874}, Berlin 1981.

\bibitem{Ed}  K. Eda, On $\Z$-kernel groups, Archiv der Mathematik
 {\bf 41}, 289 -- 293 (1983).


\bibitem{EO}  K. Eda, H. Ohta, On abelian groups of
integer-valued continuous functions, their $\Z$-dual and
$\Z$-reflexivity, in {\it Abelian Group Theory}, pp. 241 -- 257,
Gordon and Breach, London 1986.

\bibitem{EM}  P. Eklof, A. Mekler, Almost free modules,
Set-theoretic methods, North-Holland, Amsterdam 1990.


\bibitem{Fu} L. Fuchs, Infinite abelian groups - Volume 1,2
Academic Press, New York 1970, 1973.

\bibitem{GM}  R. G\"obel, W. May, Independence in completions and
endomorphism algebras, Forum Mathematicum {\bf 1}, 215 -- 226.
(1989).

\bibitem{GS} R. G\"obel, S. Shelah, Some nasty reflexive groups,
to appear in Mathematische Zeitschrift

\bibitem{GS1} R. G\"obel, S. Shelah,  Reflexive subgroups of
the Baer-Specker group and Martin's axiom, in preparation

\bibitem{J}  T. Jech, Set theory, Academic Press, New York 1978

\bibitem{S} S. Shelah, Non Structure Theory, Oxford University Press
(2000) in preparation

\end{thebibliography}
\end{document}